\documentclass[a4paper]{article}
\usepackage{amssymb}
\usepackage{amsthm}
\usepackage{amsmath}
\usepackage{latexsym}
\usepackage[T1]{fontenc}
\usepackage[latin1]{inputenc}

\def \eop {\hbox{}\nobreak\hfill \vrule width 2.0mm height 1.8mm depth 0mm
\par \goodbreak \smallskip}

\newcommand{\integ}[2]{\displaystyle \int_{#1}^{#2}}
\newcommand{\dint}{\displaystyle\int}

\begin{document}
\newtheorem{definition}{Definition}[section]
\newtheorem{theorem}{Theorem}[section]
\newtheorem{proposition}{Proposition}[section]
\newtheorem{lemma}{Lemma}[section]
\newtheorem{remark}{Remark}[section]
\newtheorem{corollary}{Corollary}[section]
\newtheorem{example}{Example}[section]
\def \ce{\centering}
\def \bop {\noindent\textbf{Proof}}
\def \eop {\hbox{}\nobreak\hfill
\vrule width 2mm height 2mm depth 0mm
\par \goodbreak \smallskip}
\def \R{\mathbb{R}}
\def \Q{\mathbb{Q}}
\def \N{\mathbb{N}}
\def \P{\mathbb{P}}
\def \E{I\!\!E}
\def \T{\mathbb{T}}
\def \H{\mathbb{H}}
\def \L{\mathbb{L}}
\def \Z{\mathbb{Z}}
\def \bf{\textbf}
\def \it{\textit}
\def \sc{\textsc}
\def \ni {\noindent}
\def \sni {\ss\ni}
\def \bni {\bigskip\ni}
\def \ss {\smallskip}
\def \F{\mathcal{F}}
\def \g{\mathcal{g}}
\def \eop {\hbox{}\nobreak\hfill
\vrule width 2mm height 2mm depth 0mm
\par \goodbreak \smallskip}
\title{Generalized Snell Envelope as a Minimal Solution of BSDE With Lower Barriers}
\author{E. H. Essaky$^1$  \quad \quad M. Hassani$^1$ \quad \quad Y. Ouknine$^2$\\\\
$^1$ Universit\'{e} Cadi Ayyad\\ Facult\'{e} Poly-disciplinaire\\
Laboratoire de Modélisation et Combinatoire\\
D\'{e}partement de Math\'{e}matiques et d'Informatique\\ B.P. 4162,
Safi, Maroc.\\ e-mails : essaky@ucam.ac.ma
\hspace{.2cm}medhassani@ucam.ac.ma \hspace{.2cm} \\ \\
$^2$ Universit\'{e} Cadi Ayyad\\ Facult\'{e} des Sciences Semlalia\\
D\'{e}partement de Math\'{e}matiques\\ B.P. 2390,  Marrakech, Maroc.\\
e-mail : ouknine@ucam.ac.ma}
\date{}
\maketitle \maketitle \footnotetext[1]{This work is supported by
Hassan II Academy of Science and technology, Action Int\'egr\'ee
MA/10/224 and Marie Curie ITN n$^{\circ}$ 213841-2.}

\begin{abstract}
The aim of this paper is to characterize the snell envelope of a
given ${\cal P}-$measurable process $l :=(l_t)_{0\leq t\leq T}$ as
the minimal solution of some backward stochastic differential
equation with lower general reflecting barriers and to prove that
this minimal solution exists.
\end{abstract}

\ni \textbf{Keys Words:} Backward stochastic differential equation;
comparison theorem; Snell envelope.
\medskip

\ni \textbf{AMS Classification}\textit{(1991)}\textbf{: }60H10,
60H20.
\bigskip

\section{Introduction and notations}
Let $(\Omega, {\cal F}, (\F_t)_{t\leq T}, P)$ be a stochastic basis
on which is defined a Brownian motion $(B_t)_{t\leq T}$ such that
$({\cal F}_t)_{t\leq T}$ is the natural filtration of $(B_t)_{t\leq
T}$ and ${\cal F}_0$ contains all $P$-null sets of $\cal F$.  Note
that $({\cal F}_t)_{t\leq T}$ satisfies the usual conditions,
\it{i.e.} it is right continuous and complete.

Let us first introduce the following notations :
\medskip

$\bullet$ $\cal P$ is the sigma algebra of ${\cal F}_t$-predictable
sets on $\Omega\times [0,T].$

$\bullet$ ${\cal D}$ is the set of ${\cal P}$-measurable and right
continuous with left limits (\it{rcll} for short) processes
$(Y_t)_{t\leq T}$ with values in $\R$.

$\bullet$ For a given process $Y\in {\cal D}$, we denote :
$Y_{t-}=\displaystyle{\lim_{s\nearrow t}Y_s}, t\leq T$
$(Y_{0^-}=Y_0)$, and $\Delta_s Y = Y_s - Y_{s-}$ the size of its
jump at time $s$.

$\bullet$ ${\cal K} := \{ K\in  {\cal D}\quad :\quad K \quad\mbox{is
nondecreasing and }\,\, K_0  =0\}$.





$\bullet$ ${\cal L}^{2,d}$ the set of $\R^d$-valued and $\cal
P$-measurable processes $(Z_t)_{t\leq T}$ such that
$$\integ{0}{T}|Z_s|^2ds<\infty, P- a.s.$$

The aim of this paper is to characterize the snell envelope of a
given ${\cal P}-$measurable process $l :=(l_t)_{0\leq t\leq T}$ as
the minimal solution of some reflected BSDE with lower barriers
(RBSDE for short). 

Let $l :=(l_t)_{0\leq t\leq T}$ be an ${\cal F}_t$-adapted right
continuous with left limits (\it{rcll} for short) process with
values in $\R$ of class ${D}[0,T]$, that is the family
$(l_{\nu})_{\nu\in \cal T}$ is uniformly integrable, where $\cal T$
is the set of all ${\cal F}_t$-stopping times $\nu$, such that
$0\leq \nu\leq T$. The Snell envelope ${\cal S}_t(l)$ of $l
:=(l_t)_{0\leq t\leq T}$ is defined as
\begin{equation}
\label{snell}
{\cal S}_{t}\left( l \right) =ess\sup_{\nu\in {\cal T}_{t}}I\!\!E%
\left[ l _\nu |{\cal F}_{t}\right],
\end{equation}
where ${\cal T}_{t}$ is the set of all stopping times valued between
$t$ and $T$. According to the work of Mertens (see \cite{dm}),
${\cal S}$ is the smallest \it{rcll}-supermartingale of class ${D
}[0,T]$ which dominates the process $l$, \it{i.e.}, $P$-a.s,
$\forall t\leq T$,\,\, $ l_t\leq {\cal S}_{t}\left( l \right)$.

Suppose now that $l$ is neither of class ${D }[0,T]$ nor a \it{rcll}
process but just ${\cal P}-$measurable, it is natural to ask whether
 we can define the smallest
local supermartingale which dominates the process $l$? In order to
give a positive answer to this question, let $L\in {\cal D}$ and
$\delta\in {\cal K}$ and assume that there exists a local martingale
$M_t = M_0 + \dint_0^t \kappa_s dB_s$ such that $P-$a.s.,
$$ L_t\leq M_t \,\,\mbox{on}\,\, [0, T[ \,\,\mbox{and} \quad l_t\leq
M_{t} \,\,d\delta_t-a.e.\,\,\mbox{on}\,\, [0,
T]\,\,\,\mbox{and}\,\,\, l_T\leq M_T. $$ Theorem \ref{thee1} states
that $Y$ the minimal solution of the following RBSDE with lower
barriers $L$ and $l$,
\begin{equation}
\label{0} \left\{
\begin{array}{ll}
(i) & 
 Y_{t}=L_T
+\integ{t}{T}dK_{s}^+ -\integ{t}{T}Z_{s}dB_{s}\,, t\leq T,
\\ (ii)&
\forall t\in[0,T[,\,\, L_t \leq Y_{t},
\\ (iii)&  \mbox{on}\,\,]0,T],\,\, l_t\leq Y_{t-},\,\,\, d\delta_t-a.e.\\
(iv)& \forall L^*\in {\cal D}\quad\mbox{satisfying}\quad  \forall t<
T, L_t\leq L_t^* \leq Y_{t}\quad\mbox{and}
\\ &\mbox{on}\,\,]0,T],\,\,
 l_t\leq L_{t-}^*,\,\,\, d\delta_t-a.e.\\ & \mbox{we have}  \integ{0}{T}(
Y_{t-}-L_{t-}^*)
dK_{t}^+=0,\,\, \mbox{a.s.}, \\
(v)& Y\in {\cal D}, \quad K^+\in {\cal K}, \quad Z\in {\cal
L}^{2,d},
\end{array}
\right. \end{equation}
is the smallest
\it{rcll} local supermartingale satisfying
$$
\forall t\in [0, T[,\,\,L_t\leq Y_t, \,\,l_t\leq Y_{t-} \,\,d\delta-
a.e., \,\,\,\mbox{on} \,\,\,[0 ,T]\,\,\,\mbox{and}\,\,\, l_T\leq
Y_T.
$$  The process $Y$ will be called later the generalized Snell envelope associated to $L,
l$ and $\delta$ and it will be denoted by $\mathcal{S}_.(L,
l,\delta, l_T)$. It is worth mentioning here that when the process
$l$ is bounded and progressively measurable and $\delta$ is the
Lebesgue measure, L. Stettner and J. Zabczyk characterize the strong
Snell envelope $V$, which is the smallest right continuous
non-negative supermartingale such that $V\geq l,\,\,$$dtdP-$a.s., as
the limit of some non-linear equation.


As by product, if we suppose that there exist $L\in {\cal D}$ and
$M$ a local martingale such that $L_{t}\leq l_t\leq M_t$,\,\
$dt-$a.e. and $l_T\leq M_T$. We prove that $Y$ the minimal solution
of the following reflected BSDE
\begin{equation}
\label{eqmin} \left\{
\begin{array}{ll}
(i) &
 Y_{t}=L_T+\integ{t}{T}dK_{s}^+ -\integ{t}{T}Z_{s}dB_{s}\,, t\leq T,
\\ (ii)& \mbox{on}\,\, ]0, T],\,\,L_t \leq Y_{t},\,\, dt- a.e
\\ (iii)&   \forall L^*\in {\cal D}\quad\mbox{satisfying}\,\, L_t\leq
L_{t}^*\leq Y_{t}\,\,\, dt- a.e.\,\, \mbox{we have}\\
&\integ{0}{T}( Y_{t-}-L_{t-}^*)
dK_{t}^+=0,\,\, \mbox{a.s.}, \\
(v)& Y\in {\cal D}, \quad K^+\in {\cal K}, \quad Z\in {\cal
L}^{2,d},
\end{array}
\right.
\end{equation}
 is the smallest \it{rcll} local supermartingale bounding the given
process $l :=(l_t)_{0\leq t\leq T}$, \it{i.e.}
$$
l_t\leq Y_t,\,\,dt- a.e\,\,\,\,\mbox{and}\,\,\,\, l_T\leq Y_T.
$$
We shall prove later that equation (\ref{0}) has a minimal solution.
We shall also
characterize the solution $Y$ as the generalized snell envelope
$\mathcal{S}_.(L) = \mathcal{S}_.(L, l,\delta, L_T)$ and we shall
show that the generalized snell enveloppe $\mathcal{S}_.(L, 0,0,
L_T)$ coincides with the usual snell envelope defined by equality
(\ref{snell}) if the process $L$
is of class $\it D[0,T]$.\\

We need also the following notations :\\

\ni$\bullet$ For a set $B$, we denote by $B^c$ the
complement of $B$ and $1_B$ denotes the indicator of $B$.\\
$\bullet$ For each $(a, b)\in\R^2$, $a\wedge b = \min(a, b)$\,\, and
\,\,$a\vee b = \max(a, b)$.

\ni$\bullet$ For all $(a, b, c)\in\R^3$ such that $a\leq c$,\,\,
$a\vee b\wedge c = \min(\max(a,b), c) = \max(a, \min(c,b))$.\\

Throughout the paper we introduce the following data :\\

\ni$\bullet$ $\xi$ is an ${\cal F}_T$-measurable one dimensional
random
variable.\\

\ni$\bullet$ $L:=\left\{ L_{t},\,0\leq t\leq T\right\}$ is a barrier
which belongs to ${\cal D}$.
\\

\ni$\bullet$ $l:=\left\{ l_{t},\,0\leq t\leq T\right\}$ is a $ {\cal
P}-${measurable} process.

\ni$\bullet$ $\delta\in {\cal K}$.

\ni$\bullet$ ${\cal M} ={\cal M}(L, l,\delta, \xi)$ is the set of
\it{rcll} local supermartingale $V_t = V_0 -A_t +\dint_0^t \chi_s
dB_s$, where $A\in {\cal K}$ and $\chi \in {\cal L}^{2, d}$ such
that
$$
L_t\leq V_t, \,\,l_t\leq V_{t-} \,\,d\delta_t-
a.e.\,\,\,\mbox{and}\,\,\, \xi\leq V_T.
$$
We should note here that if $V_t = V_0 -A_t +\dint_0^t \chi_s
dB_s\in {\cal M}$, then we have
\begin{enumerate}
\item $V_t +1\in {\cal M}.$
\item $V_t+A_t = V_0 +\dint_0^t \chi_s
dB_s\in {\cal M}.$
\end{enumerate}

\section{Preliminaries}
In view of clarifying this issue, we recall some results concerning
generalized reflected BSDEs (GRBSDE for short) with two \it{rcll}
obstacles. We present both the existence and comparison theorem for
minimal solutions of this kind of equations. Those results will play
a crucial role in our proofs (see \cite{EHO} for more details). We
should note here that the notion of reflected BSDE with two
obstacles has been first introduced by Civitanic and Karatzsas
\cite{CK}.

\subsection{Existence of a minimal solutions for GRBSDE}
Let us recall first the following definition of two singular
measures.
\begin{definition}
Let $K^1$ and $K^2$ be two processes in ${\cal K}$. We say that :\\
\ni $K^1$ and $K^2$ are singular if and only if there exists a set
$D\in {\cal P}$ such that
$$
\E\dint_0^T 1_D(s,\omega) dK^1_s(\omega) = \E\dint_0^T
1_{D^c}(s,\omega) dK^2_s(\omega) =0.
$$
 This is denoted by $dK^1 \perp dK^2$.
\end{definition}
Let us now define the notion of solution of the GRBSDE with two
obstacles $L$ and $U$. For this reason, let : \\
$\bullet$ $g : [ 0,T]\times \Omega \times \R\longrightarrow \R$ be a
 function such that
$$
\forall y\in \R,\,\,\, (t,\omega)\longmapsto g(t, \omega,
L_{t-}(\omega)\vee y \wedge U_{t-}(\omega))\,\,\mbox{ is }\,\, {\cal
P}-\mbox{measurable}.
$$
$\bullet$ $U:=\left\{ U_{t},\,0\leq t\leq T\right\}$ be a barrier
 such that $L_t\leq U_t$, $\forall t\in
 [0,T[$. 
\begin{definition}\label{def01}
\begin{enumerate}
 \item We say that $(Y,Z,K^+,K^-):=( Y_{t},Z_{t},K_{t}^+,K_{t}^-)_{t\leq T}$
is a solution of the generalized reflected BSDE, associated with the
data $(\xi, g, \delta, L, U)$, if the following hold :
\begin{equation}
\label{eq001} \left\{
\begin{array}{ll}
(i) & 
 Y_{t}=\xi
+\dint_t^Tg(s, Y_{s-})d\delta_s  +\integ{t}{T}dK_{s}^+
-\integ{t}{T}dK_{s}^- -\integ{t}{T}Z_{s}dB_{s}\,, t\leq T,
\\ (ii)&
\forall t\in[0,T[,\,\, L_t \leq Y_{t}\leq U_{t},\\  (iii)&
\integ{0}{T}( Y_{t-}-L_{t-})
dK_{t}^+= \integ{0}{T}( U_{t-}-Y_{t-}) dK_{t}^-=0,\,\, \mbox{a.s.}, \\
(iv)& Y\in {\cal D}, \quad K^+, K^-\in {\cal K}, \quad Z\in {\cal
L}^{2,d},  \\ (v)& dK^+\perp  dK^-.
\end{array}
\right. \end{equation}
\item We say that the GRBSDE (\ref{eq001}) has a minimal
solution $(Y_t ,Z_t ,K^+_t , K_t^- )_{t\leq T}$ if for any other
solution $(Y_t^{'} ,Z_t^{'} ,K'^{+}_t , K'^{-}_t )_{t\leq T}$ of
(\ref{eq001}) we have for all $t \leq T$, $Y_t\leq Y_t^{'}$,
$P$-a.s.

\end{enumerate}
\end{definition}

We introduce also the following assumption : \\ \\ \ni $(\bf{H})$
The function $g$ and the barrier $U$ satisfy the following :
\begin{itemize}
\item[\bf{(a)}]There exists $\beta \in L^0(\Omega, L^1([0,T], \delta(dt),
\R_+))$ such that : $  \forall y\in\R,\quad |g(t, \omega,
L_{t-}(\omega)\vee y \wedge U_{t-}(\omega))| \leq
\beta_t(\omega),\quad \delta(dt)P(d\omega)-$a.e.
\item[\bf{(b)}] $\delta(dt)P(d\omega)-$a.e.,\,\,the function\,\,
$y\longmapsto g(t, \omega, L_{t-}(\omega)\vee y \wedge
U_{t-}(\omega))$ is continuous.
\item[\bf{(c)}] The barrier $U$ is a \it{rcll} local supermartingale, \it{i.e.} there exist $\alpha\in {\cal K}$ and $\gamma\in {\cal L}^{2,d}$
such that $U_t =U_0 -\alpha_t+\dint_0^t \gamma_sdB_s$.

\end{itemize}\vspace{0.2cm}
The following theorem has already been proved in \cite{EHO}. We
should note here that the barriers $L$ and $U$ are \it{rcll}, the
continuous case has been studied in \cite{EH}.
\begin{theorem}\label{thee1}
If assumption $(\bf{H})$ holds then the GRBSDE (\ref{eq001}) has a
minimal solution.
\end{theorem}
\subsection{Comparison theorem for minimal solutions } Let us now
recall the following comparison theorem which plays a crucial rule
in the proof of the existence of solutions for RBSDE. The proof of
this comparison theorem is based on an exponential change and an
approximation scheme, see \cite{EHO}. Let $(Y, Z, K^+, K^-)$ be the
minimal solution for the following GRBSDE
\begin{equation}
\label{eq0q} \left\{
\begin{array}{ll}
(i) & 
 Y_{t}=\xi+\dint_t^Tg(s, Y_{s-})d\delta_s
+\integ{t}{T}dK_{s}^+ -\integ{t}{T}dK_{s}^-
-\integ{t}{T}Z_{s}dB_{s}\,, t\leq T,
\\ (ii)&
\forall t\in[0,T[,\,\, L_t \leq Y_{t}\leq U_{t},\\  (iii)&
\integ{0}{T}( Y_{t-}-L_{t-})
dK_{t}^+= \integ{0}{T}( U_{t-}-Y_{t-}) dK_{t}^-=0,\,\, \mbox{a.s.}, \\
(iv)& Y\in {\cal D}, \quad K^+, K^-\in {\cal K}, \quad Z\in {\cal
L}^{2,d},  \\ (v)& dK^+\perp  dK^-.
\end{array}
\right. \end{equation} Let $(Y', Z', K'^+, K'^-)$ be a solution for
the following GRBSDE
\begin{equation}
\label{eq1q} \left\{
\begin{array}{ll}
(i) & 
 Y'_{t}=\xi'+\dint_t^TdA'_s+\integ{t}{T}dK'^{+}_{s} -\integ{t}{T}dK'^{-}_{s}
-\integ{t}{T}Z'_{s}dB_{s}\,, t\leq T,
\\ (ii)&
\forall t\in[0,T[,\,\, L'_t \leq Y'_{t}\leq U'_{t},\\  (iii)&
\integ{0}{T}( Y'_{t-}-L'_{t-})
dK'^{+}_{t}= \integ{0}{T}( U'_{t-}-Y'_{t-}) dK'^{-}_{t}=0,\,\, \mbox{a.s.}, \\
(iv)& Y'\in {\cal D}, \quad K'^+, K'^-\in {\cal K}, \quad Z'\in
{\cal L}^{2,d},  \\ (v)& dK'^+\perp  dK'^-,

\end{array}
\right. \end{equation} where $A'$ is a process in ${\cal K}$, $L'$
and $U'$ are two barriers
 which belong to ${\cal D}$ such that $L'_t\leq U'_t$, $\forall t\in
 [0,T[$. \\ \ni Assume moreover that for every $t\in[0, T]$
\begin{itemize}
\item[(a)] $\xi\leq \xi'$.
\item[(b)] $Y'_t\leq U_t$,\,\, $L'_t\leq Y_t$, $\forall t\in
 [0,T[$.
 \item[(c)] $g(s, Y'_{s-})d\delta_s\leq dA'_s$ on $[0, T]$.
\end{itemize}
\begin{theorem} (Comparison theorem for minimal solutions, see \cite{EHO})
\label{th2} Assume that the above assumptions hold then we have :
\begin{enumerate}
\item $Y_t\leq Y'_t$, for every $t\in [0,T]$, $P-$a.s.
\item $
1_{\{U'_{t-} = U_{t-}\}}dK^{-}_t \leq dK'^{-}_t\,\,\,
\mbox{and}\,\,\, 1_{\{L'_{t-} = L_{t-}\}}dK'^{+}_t \leq dK^{+}_t. $
\end{enumerate}
\end{theorem}
\section{Generalized Snell envelope as a solution of some RBSDE}
In this section, we prove an existence result of a minimal solution
for some reflected BSDE with lower barriers. We shall also
characterize this minimal solution $Y$ as the smallest \it{rcll}
local supermartingale satisfying
$$
\forall t\in [0, T[,\,\,L_t\leq Y_t, \,\,l_t\leq Y_{t-}
\,\,d\delta_t- a.e., \,\,\,\mbox{on} \,\,\,[0
,T]\,\,\,\mbox{and}\,\,\, \xi\leq Y_T.
$$
\ni Let us now introduce the definition of our RBSDE with lower
obstacles.
\begin{definition}\label{def1}
\begin{enumerate}
 \item We call $(Y,Z,K^+):=( Y_{t},Z_{t},K_{t}^+)_{t\leq T}$
a solution of the RBSDE, associated with the data $(\xi,L,
l,\delta)$, if the following hold :
\begin{equation}
\label{eq0} \left\{
\begin{array}{ll}
(i) & 
 Y_{t}=\xi
+\integ{t}{T}dK_{s}^+ -\integ{t}{T}Z_{s}dB_{s}\,, t\leq T,
\\ (ii)&
\forall t\in[0,T[,\,\, L_t \leq Y_{t},
\\ (iii)&  \mbox{on}\,\,]0,T],\,\, l_t\leq Y_{t-},\,\,\, d\delta_t-a.e.\\
(iv)& \forall L^*\in {\cal D}\quad\mbox{satisfying}\quad  \forall t<
T, L_t\leq L_t^* \leq Y_{t} \quad\mbox{and}
\\ &\mbox{on}\,\,]0,T],\,\,\ l_t\leq
L_{t-}^*,\,\,\, d\delta_t-a.e.\\ & \mbox{we have}  \integ{0}{T}(
Y_{t-}-L_{t-}^*)
dK_{t}^+=0,\,\, \mbox{a.s.}, \\
(v)& Y\in {\cal D}, \quad K^+\in {\cal K}, \quad Z\in {\cal
L}^{2,d}.
\end{array}
\right. \end{equation}
\item We say that the RBSDE (\ref{eq0}) has a minimal
solution $(Y_t ,Z_t ,K^+_t )_{t\leq T}$ if for any other solution
$(Y_t^{'} ,Z_t^{'} ,K^{'+}_t)_{t\leq T}$ of (\ref{eq0}) we have for
all $t \leq T$, $Y_t\leq Y_t^{'}$, $P$-a.s.
\end{enumerate}
\end{definition}
\subsection{Main result}
Let $L\in {\cal D}$, $\xi\in L^0(\Omega)$, $l\in L^0(\Omega\times
[0, T])$ and $\delta\in {\cal K}$. We assume the following
hypothesis :\\

\ni \bf{(A)} There exists a local martingale $M_t = M_0 + \dint_0^t
\kappa_s dB_s$ such that $P-$a.s., $ L_t\leq M_t \,\,\mbox{on}\,\,
[0, T[ \,\,\mbox{and} \\ l_t\leq M_{t}
\,\,d\delta_t-a.e.\,\,\mbox{on}\,\, [0, T]\,\,\,\mbox{and}\,\,\,
\xi\leq M_T. $ This is equivalent to ${\cal M} \neq \emptyset$.\\

The main result of this paper is the following.
\begin{theorem}\label{thee1}
If assumption $(\bf{A})$ hold then the RBSDE (\ref{eq0}) has a
minimal solution $(Y_t ,Z_t ,K^+_t )_{t\leq T}$. Moreover $Y$ is the
smallest \it{rcll} local supermartingale satisfying
$$
\forall t\in [0, T[,\,\,L_t\leq Y_t, \,\,l_t\leq Y_{t-}
\,\,d\delta_t- a.e., \,\,\,\mbox{on} \,\,\,[0
,T]\,\,\,\mbox{and}\,\,\, \xi\leq Y_T.
$$
We say that $Y$ is the generalized Snell envelope associated to $L,
l, \delta$ and $\xi$. We denote it by $\mathcal{S}_.(L, l,\delta,
\xi)$.
\end{theorem}
\subsubsection{Auxiliary penalized equation}
 \ni Let $V_t = V_0 -A_t +\dint_0^t \chi_s
dB_s\in {\cal M}$.
 Let also
$({Y}^{(n,V)},{Z}^{(n,V)},{K}^{(n,V)+},{K}^{(n,V)-})$ be the minimal
solution of the following penalized RBSDE with two \it{rcll}
barriers
\begin{equation}
\label{eq00} \left\{
\begin{array}{ll}
(i) &
 {Y}^{(n,V)}_{t}=\xi + n\dint_t^T(l_s
-{Y}^{(n,V)}_{s-})^+d\delta_s+\integ{t}{T}d{K}_{s}^{(n,V)+}
\\
&\qquad\quad-\integ{t}{T}d{K}_{s}^{(n,V)-}
-\integ{t}{T}{Z}_{s}^{(n,V)}dB_{s}\,, t\leq T,
\\ (ii)&
\forall t\in[0,T[,\,\, L_t \leq {Y}^{(n,V)}_{t}\leq V_t,\\
(iii)& \integ{0}{T}( {Y}^{(n,V)}_{t-}-L_{t-})
d{K}_{t}^{{(n,V)}+}= \integ{0}{T}( V_{t-}-{Y}^{(n,V)}_{t-}) d{K}_{t}^{{(n,V)}-}=0,\,\, P-\mbox{a.s.}, \\
(iv)& {Y}^{(n,V)}\in {\cal D}, \quad {K}^{{(n,V)}+},
{K}^{{(n,V)}-}\in {\cal K}, \quad {Z}^{(n,V)}\in {\cal L}^{2,d},
\\ (v)& d{K}^{{(n,V)}+}\perp  d{K}^{{(n,V)}-}.
\end{array}
\right. \end{equation} 
We should mention here that the minimal solution to (\ref{eq00})
exists according to Theorem
\ref{thee1} (see \cite{EHO} for the proof).\\\\
 \ni Our objective now is to prove that ${Y}^{(n,V)}$ does not
depend on $V\in {\cal M}$ and converges to some ${Y} $ which belongs
to $\in {\cal M}$. This means that the process $Y$ is the smallest
\it{rcll} local supermartingale satisfying
$$
\forall t\in [0, T[,\,\,L_t\leq Y_t, \,\,l_t\leq Y_{t-}
\,\,d\delta_t- a.e., \,\,\,\mbox{on} \,\,\,[0
,T]\,\,\,\mbox{and}\,\,\, \xi\leq Y_T.
$$
It follows from comparison theorem \ref{th2}, applied to
${{Y}}^{(n,V)}$ and $V_t$ (we can also apply Tanaka's formula to the
process $(V_t-{Y}^{(n,V)}_t)^+ = (V_t-{Y}^{(n,V)}_t)$), that for
every $n\in\N$ $d{K}^{{(n,V)}-}= 0$. Hence
$({Y}^{(n,V)},{Z}^{(n,V)},{K}^{(n,V)+})$ is the minimal solution of
the following GRBSDE
\begin{equation}
\label{eq000} \left\{
\begin{array}{ll}
(i) &
 {Y}^{(n,V)}_{t}=\xi + n\dint_t^T(l_s
-{Y}^{(n,V)}_{s-})^+d\delta_s+\integ{t}{T}d{K}_{s}^{(n,V)+}
\\
&\qquad\quad -\integ{t}{T}{Z}_{s}^{(n,V)}dB_{s}\,, t\leq T,
\\ (ii)&
\forall t\in[0,T[,\,\, L_t \leq {Y}^{(n,V)}_{t},\\
(iii)& \integ{0}{T}( {Y}^{(n,V)}_{t-}-L_{t-})
d{K}_{t}^{{(n,V)}+}=0,\,\, P-\mbox{a.s.}, \\
(iv)& {Y}^{(n,V)}\in {\cal D}, \quad {K}^{{(n,V)}+}\in {\cal K},
\quad {Z}^{(n,V)}\in {\cal L}^{2,d}.
\end{array}
\right. \end{equation} Moreover, for every $V\in {\cal M}$ and all
$(n,t)\in\N\times
[0,T]$, ${Y}^{(n,V)}_{t}\leq V_t$.\\
Since ${Y}^{(n,M)}$ is also the minimal solution of $(\ref{eq000})$,
then for every $V$, ${Y}^{(n,V)}={Y}^ {(n,M)}$. From now on we
denote the solution of $(\ref{eq000})$
by $({Y}^{n},{Z}^{n},{K}^{n+})$. \\
Now by using comparison theorem \ref{th2} we get, for every $V\in
{\cal M}$, that
\begin{equation}\label{equ1}
 L_t \leq {Y}^{n}_{t}\leq {Y}^{{n+1}}_{t}\leq
 V_t.
\end{equation}
Now let us set \begin{equation}\label{equ10}
\begin{array}{ll}
& 
{Y}_{t}=\displaystyle{\sup_{n}}{Y}^{n}_{t} \,\,\mbox{and}\,\,
{Y}_{t}^{\,-}=\displaystyle{\sup_{n}}{Y}^{n}_{t-}.
\end{array}
\end{equation}
The following results guarantee that the process $Y$ is the smallest
\it{rcll} local supermartingale satisfying
$$
\forall t\in [0, T[,\,\,L_t\leq Y_t, \,\,l_t\leq Y_{t-}
\,\,d\delta_t- a.e., \,\,\,\mbox{on} \,\,\,[0
,T]\,\,\,\mbox{and}\,\,\, \xi\leq Y_T.
$$
By letting $n$
to infinity in (\ref{equ1}) and using assumption $(\bf{A})$ we have
the following.
\begin{lemma}\label{lem0}
 For every $t\in [0, T]$ we have for every $V\in {\cal M}$,
$$
L_t \leq {Y}_{t}\leq V_t \quad \mbox{on}\quad [0,T[ \quad
\mbox{and}\quad L_{t-} \leq {Y}_{t}^{-}\leq V_{t-}\quad
\mbox{on}\quad ]0,T].
$$
\end{lemma}
\begin{proposition}\label{pro1} The process ${Y}$ defined by (\ref{equ10}) satisfy the
following properties :
\begin{enumerate}
\item ${Y}$ is a \it{rcll} local supermartingale and ${Y}_{t}^-\leq {Y}_{t-}$, for every $t\in ]0,T]$.
\item $l_t\leq {Y}_{t}^-, \,\,\, d\delta_t-a.e.,$ \mbox{on}
$]0,T]$.
\end{enumerate}
\ni In particular it follows that $Y$ belongs to ${\cal M}$.
\end{proposition}
\bop. \it{1.} Recall that $M_t = M_0 + \dint_0^t \kappa_s dB_s \in
{\cal M}$. We have
$$
\begin{array}{ll}
& {Y}_{t}^n-M_t\\ & =\xi- M_T + n\dint_t^T
(l_s-{Y}_{s-}^n)^+d\delta_s + \dint_t^T d{K}_{s}^{n+}+
\dint_t^T(\overline{Z}_{s}^n-\kappa_s)dB_s.
\end{array}
$$

\ni Let $(\tau_i)_{i\geq 1}$ be the family of stopping times defined
by
\begin{equation}\label{tau}
\tau_i =\inf\{s\geq 0 : M_{s} -L_s\geq i+M_0- L_0\}\wedge T.
\end{equation}
Note that $\tau_i >0,\,\,P-$a.s., for every $i\geq 1$.
By using a localization procedure we have for every $ i\geq 1$ and $n\geq 0$
\begin{equation}\label{equ6}
\begin{array}{ll}
  & \E(M_0-{Y}_{0}^n)+n
  \E\dint_0^{{\tau_{i}-}}(l_s-{Y}_{s-}^n)^+d\delta_s+\E
{K}_{\tau_i-}^{n+} \leq i+\E(M_0- L_0).
\end{array}
\end{equation}
Put
\begin{equation}\label{equ7}
\begin{array}{ll}
  & \mathcal{M}_t^n ={Y}_{t}^n-M_t,
\\ & ^i\mathcal{M}_t^n = \mathcal{M}_t^n 1_{\{t<\tau_i\}} + \mathcal{M}_{{\tau_i}-}^n
1_{\{t\geq\tau_i\}},
\end{array}
\end{equation}
we have
$$
-i-\E(M_0- L_0)\leq\, ^i\mathcal{M}_t^n\leq 0\quad\mbox{and}\quad
^i\mathcal{M}_t^n\leq ^i\mathcal{M}_t^{n+1} \quad\mbox{and}\quad
t\rightarrow\, ^i\mathcal{M}_t^{n}\,\,\mbox{is a \it{rcll}
supermartingale}.
$$
It follows then from Dellacherie and Meyer \cite{dm} that
$\displaystyle{\sup_{n}}\,\, ^i\mathcal{M}_t^{n}$ is also a
\it{rcll} process supermartingale process. Since
$P\bigg[\displaystyle{\bigcup_{i=1}^{\infty}(\tau_i= T)}\bigg] =1,$
it follows that
${Y}_t$ is a \it{rcll} local supermartingale on $[0, T]$.\\ 
Now since for every $s\in ]0, T]$ and $n\in\N$,\,\,${Y}_{s-}^n\leq
{Y}_{s-}$, it follows that $
{Y}_{s}^{\,\,-}\leq {Y}_{s-}$.\\\\
\ni \it{2.} On another hand, by letting $n$ to infinity in
inequality (\ref{equ6}) and using Fatou's lemma it follows that
$$
\E\dint_0^{{\tau_{i}}-} (l_s-{Y}_{s}^-)^+d\delta_s =0.
$$
Hence
$$
(l_s-{Y}_{s}^-)^+ =0\,\, d\delta_s- a.e. \,\,\mbox{on}\,\, [0, T[.
$$
Assume now that ${Y}_{T}^- < l_T$ and $\Delta_T\delta >0$. It
follows from \cite{EHO}, that for every $V\in {\cal M}$
$$
{Y}_{T-}^n = L_{T-}\vee [\xi +n(l_T-{Y}_{T-}^n
)^+\Delta_T\delta]\wedge V_{T-}\geq [\xi +n(l_T-{Y}_{T-}
)^+\Delta_T\delta]\wedge V_{T-}.
$$
We get ${Y}_{T}^{\,\,-} = V_{T}$, which is absurd since $V_t +1\in
{\cal M}$. Consequently
$$
l_s\leq{Y}_{s}^-\,\,\, d\delta_s- a.e. \,\,\mbox{on}\,\, [0, T].
$$
 The proof of Proposition \ref{pro1} is finished. \eop

\subsubsection{Proof of the main result}
\bop\,\,\bf {of Theorem \ref{thee1}}. Let $L^*\in {\cal D}$ be such
that for every $t\in [0,T]$,\, $L_t\leq L^*_t\leq Y_t$ and $l_t\leq
L_{t-}^*\,\,d\delta_t- a.e.$. Let also $(Y^*, Z, K^{+}, K^{-})$,
which is exists according to Theorem \ref{thee1}, the minimal
solution of the following RBSDE
\begin{equation}
\label{eq01000} \left\{
\begin{array}{ll}
(i) & 
 Y_{t}^*=\xi
+\integ{t}{T}dK_{s}^{+} -\integ{t}{T}dK_{s}^{-}
-\integ{t}{T}Z_{s}dB_{s}\,, t\leq T,
\\ (ii)&
\forall t\in[0,T[,\,\, L_{t}^* \leq Y_{t}^*\leq Y_{t},
\\ (iii)&   \integ{0}{T}(
Y_{t-}-Y_{t-}^*)
dK_{t}^{-}= \integ{0}{T}(Y_{t-}^*-L_{t-}^*) dK_{t}^{+}=0,\,\, \mbox{a.s.}, \\
(v)& Y^*\in {\cal D}, \quad K^{+}, K^{-}\in {\cal K}, \quad Z\in
{\cal L}^{2,d},  \\ (vi)& dK^{+}\perp  dK^{-}.
\end{array}
\right. \end{equation} By the same argument as before with $V =Y$
($Y$ is the process defined in the previous subsection), one can see
that $dK^{-} =0$, hence $Y^*\in {\cal M}$. By Lemma $\ref{lem0}$ and
$(ii)$ of Equation (\ref{eq01000}) we get

$$
Y^*_s= Y_s.
$$
Henceforth
$$
(Y_{t-}-L_{t-}^*) dK_{t}^{+} = 0.
$$
Consequently, for every $V\in {\cal M}$, $(Y, Z, K^+)$ is a solution
of (\ref{eq0}). Moreover the process $Y$ is the smallest \it{rcll}
local supermartingale satisfying
$$
\forall t\in [0, T[,\,\,L_t\leq Y_t, \,\,l_t\leq Y_{t-}
\,\,d\delta_t- a.e., \,\,\,\mbox{on} \,\,\,[0
,T]\,\,\,\mbox{and}\,\,\, \xi\leq Y_T.
$$\eop

As by product we obtain the following theorem.
\begin{theorem} Let $(T_i)_{i\geq 1}$ be a sequence of stopping times such that $[|T_i|]\cap[|T_j|] =
\emptyset,\,\, \forall i\neq j$ and $\bigcup_{i\geq 1}[|T_i|] =\{(t,
\omega)\in ]0,T]\times\Omega : \Delta_t\delta(\omega)
>0\}$. Under assumption \bf{(A)}, $Y$ the minimal solution of (\ref{eq0}) is the smallest
\it{rcll} local supermartingale satisfying $P-$a.s.
$$
\forall t\in [0, T[,\,\,L_t\leq Y_t, \,\,l_t\leq Y_{t}
\,\,d\delta_t^c- a.e., \,\,\,\mbox{on} \,\,\,[0 ,T],\,\,\,\forall
i\geq 1,\,\, l_{T_i}\leq Y_{T_{i}-} \,\,\,\mbox{and}\,\,\, \xi\leq
Y_T.
$$
\end{theorem}
\subsubsection{Some properties of the generalized Snell
envelope} 
The generalized Snell envelope $Y=\mathcal{S}_.(L, l,\delta, \xi)$
solution of RBSDE (\ref{eq0}) has the following properties whose
proofs are immediate.

\begin{corollary}
\begin{enumerate}
\item $\mathcal{S}_.(L, l,\delta, \xi) = \mathcal{S}_.(L, \overline{l},\delta, \xi)$,
with $\overline{l}_s = l_s \vee L_{s-}$.
\item If $ L'\leq L$, $d\delta'\ll d\delta$,  $l'\leq l,\,\,\, d\delta'$ a.e., $\xi'\leq
\xi$ then $(L', l', \delta', \xi')$ satisfies condition \bf{(A)} and
$\mathcal{S}_.(L', l',\delta', \xi')\leq \mathcal{S}_.(L, l,\delta,
\xi).$
\item $\mathcal{S}_.(L, l,\delta, \xi)\geq  \mathcal{S}(L^{\xi})$ (with equality if $l_t\leq L_{t-}
\,\,d\delta_t- a.e., \,\,\,\mbox{on} \,\,\,[0 ,T]$) where
$\mathcal{S}(L^{\xi}) = \mathcal{S}_.(L, 0,0, \xi)$ and $L^{\xi}_t =
L_t 1_{\{t<T\}} +\xi 1_{\{t=T\}}$.
\item Put $Y=\mathcal{S}_.(L, l,\delta, \xi)$. If
$$
l_t\leq l'_t\leq Y_{t-},\,\,\,d\delta - a.e., \,\, \mbox{on}\,\, [0,
T]\,\, \mbox{and}\,\, L_t\leq L'_t\leq Y_t,\, \forall t\in [0, T[,
\,\, \mbox{and}\,\, d\delta\sim d\delta',
$$
then $\mathcal{S}_.(L, l,\delta, \xi) = \mathcal{S}_.(L',
l',\delta', \xi)$.\\

\ni In particular for every $L^*\in {\cal D}$ such that $P-$a.s., $$
L_t\leq L^*_t\leq Y_t,\, \forall t\in [0, T[,  \,\, \mbox{and}\,\,
l_t\leq L^*_{t-}\leq Y_{t-},\,\,\,d\delta_t- a.e., \,\,
\mbox{on}\,\, [0, T] \,\, \mbox{and}\,\, L^*_T =\xi
$$
we have $\mathcal{S}_.(L, l,\delta, \xi) = \mathcal{S}_.(L^*)$.
\end{enumerate}
\end{corollary}
\begin{remark} We know that if $L$ is of class $\it D$ then $L$ satisfies assumption
\bf{(A)} (see Dellacherie-Meyer \cite{dm}). In this case our
generalized snell enveloppe $\mathcal{S}_.(L) = \mathcal{S}_.(L,
0,0, L_T)$ coincides with the usual snell enveloppe $esssup_{\tau
\in {\cal T}_{t}} \E[L_{\tau}|{\cal F}_t]$, where ${\cal T}_{t}$ is
the set of all stopping times valued between $t$ and $T$, as
presented in Dellacherie-Meyer \cite{dm} and studied by several
authors.
\end{remark}
\begin{example} \ni If $\delta_t =t$ and there exist $L\in {\cal D}$ and $M$ a local martingale such that
$L_{t}\leq l_t\leq M_t$ and $\xi\leq M_T$. Let $(Y, Z, K^+)$ be the
minimal solution
 of the following RBSDE
$$
\label{eq010000} \left\{
\begin{array}{ll}
(i) & 
 Y_{t}=\xi
+\integ{t}{T}dK_{s}^+ -\integ{t}{T}Z_{s}dB_{s}\,, t\leq T,
\\ (ii)& \mbox{on}\,\, ]0, T],\,\,l_t \leq Y_{t},\,\, dt- a.e
\\ (iii)&   \forall L^*\in {\cal D}\quad\mbox{satisfying}\,\, l_t\leq
L_{t}^*\leq Y_{t}\,\,\, dt- a.e.\,\, \mbox{we have}\\
&\integ{0}{T}( Y_{t-}-L_{t-}^*)
dK_{t}^+=0,\,\, \mbox{a.s.}, \\
(v)& Y\in {\cal D}, \quad K^+\in {\cal K}, \quad Z\in {\cal
L}^{2,d},
\end{array}
\right. $$ Then $Y$ is the smallest local supermartingale such that
$$
l_t\leq Y_t,\,\,dt- a.e\,\,\,\,\mbox{and}\,\,\,\, \xi\leq Y_T.
$$
\end{example}

\end{document}